\documentclass[noautosecdot,aop,seceqn,secthm,dvips]{arximspdf}
\usepackage{mathbh}
\usepackage{graphicx}

\doi{10.1214/009117905000000477}
\volume{33}
\issue{6}
\pubyear{2005}
\firstpage{2127}
\lastpage{2148}

\makeatletter
\renewcommand{\aa}{p_n}
\renewcommand{\AA}{A}
\renewcommand{\Im}{\operatorname{Im}}
\newcommand{\noa}{\noalign{}}
\newcommand{\ds}{\displaystyle}
\newcommand{\up}{\textup}
\newtheorem{corollary}[thm]{Corollary}
\newtheorem{lemma}[thm]{Lemma}
\newtheorem{proposition}[thm]{Proposition}
\newproclaim{rem}[thm]{Remark}
\newcommand{\eref}[1]{(\ref{#1})}
\newcommand{\R}{\mathbb{R}}
\newcommand{\C}{\mathbb{C}}
\newcommand{\N}{\mathbb{N}}
\renewcommand{\H}{\mathbb{H}}
\newcommand{\U}{\mathbb{U}}
\newcommand{\diam}{\operatorname{diam}}
\newcommand{\area}{\operatorname{area}}
\newcommand{\dist}{\operatorname{dist}}
\newcommand{\eps}{\varepsilon}
\renewcommand{\P}{\mathbf{P}}
\newcommand{\md}{\vert}
\newcommand{\p}{\partial}
\newcommand{\E}{\mathbf{E}}
\newcommand{\closure}{\overline}
\newcommand{\ev}[1]{\mathcal{#1}}
\newcommand{\capacity}{\operatorname{cap}_\infty}
\newcommand{\hg}{\gamma^\phi}
\newcommand{\tg}{\tilde{\gamma}}
\newcommand{\inr}[1]{\mathrm{rad}_{#1}}
\newcommand{\Doms}{\mathfrak{D}}

\newcommand{\oo}{\bar{v}}
\newcommand{\bo}{\oo_\mathrm{end}}
\newcommand{\dhaus}{\mathrm{d}_H}
\newcommand{\bb}{p_0}
\newcommand{\rr}{\rho}
\newcommand{\cc}{\hat{c}}
\newcommand{\DX}{\nu}
\newcommand{\inE}{\overrightarrow{E}}
\newcommand{\outE}{\overleftarrow{E}}
\newcommand{\rev}{\mathrm{rev}}
\newcommand{\HH}{\mathfrak{h}}
\newcommand{\DD}{\tilde{D}}
\makeatother

\begin{document}
\begin{frontmatter}

\title{The harmonic explorer and its convergence to SLE(4)}
      
\runtitle{Harmonic explorer and $\mathrm{SLE}_4$}

\begin{aug}
\author[A]{\fnms{Oded} \snm{Schramm}\ead[label=e1]{schramm@microsoft.com}} and
\author[B]{\fnms{Scott} \snm{Sheffield}\corref{}\ead[label=e2]{sheff@math.berkeley.edu}}
\runauthor{O. Schramm and S. Sheffield}
\affiliation{Microsoft Research and
University of California, Berkeley}
\address[A]{Microsoft Research\\
One Microsoft Way\\
Redmond, Washington 98052\\
USA\\
\printead{e1}} %adresu isvedimo komanda gale!
\address[B]{Department of Mathematics\\
University of California, Berkeley\\
929 Euclid Avenue\\
Berkeley, California 98007\\
USA\\
\printead{e2}}
\end{aug}

% HISTORY:
\received{\smonth{10} \syear{2003}}
\revised{\smonth{2} \syear{2005}}

% ABSTRACT
%
\begin{abstract}
The harmonic explorer is a random grid path. Very roughly, at
each step the
harmonic explorer takes a turn to the right with
probability equal to the discrete harmonic measure
of the left-hand side of the path from a point near the end
of the current path. We prove that the
harmonic explorer converges
to SLE(4) as the grid
gets finer.
\end{abstract}

% KEYWORDS
%
%
% \sep

\end{frontmatter}

%spell_from *************** Text entry area ******************%

%s1 ###
\section{Introduction.}

\begin{figure}

\includegraphics{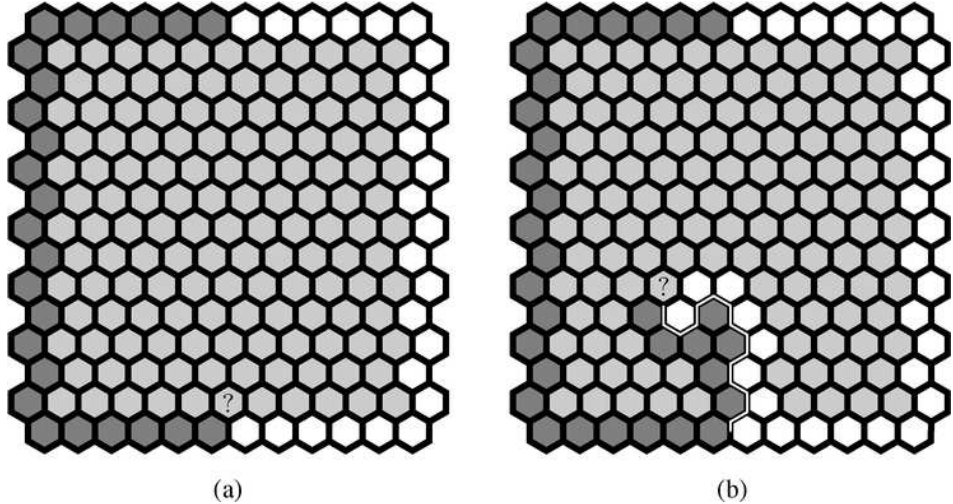}

\caption{\up{(a)} Initial left boundary faces (black), right
boundary faces (white) and
undetermined interior faces (gray). \up{(b)} A possible HE configuration
several steps later.}
\label{hesetuphexagons}
\end{figure}

Let $D$ be a simply connected subset of the hexagonal faces in the
planar honeycomb lattice. Two faces of $D$ are considered adjacent if
they share an edge. Suppose further that the boundary faces of $D$ are
partitioned into a ``left boundary'' component, colored black, and a
``right boundary'' component, colored white, in such a way that the set
of interior faces remains simply connected. (See Figure
\ref{hesetuphexagons}.) Given any black--white coloring of the faces of
$D$, there will be a~unique interface $\gamma$ separating the cluster
of black hexagons containing the left boundary from the cluster of
white hexagons containing the right boundary.

If the colors are chosen via independent Bernoulli percolation, we may
view~$\gamma$ as being
generated dynamically as follows: simply begin the path $\gamma$ at an
edge separating the left and
right boundary components; when $\gamma$ hits a black hexagon, it turns
right, and when it hits a
white hexagon, it turns left. Each time it hits a hexagon whose color
has yet to be determined, we
choose that hexagon's color with a coin toss.

The {\it harmonic explorer} (HE)
is a random interface generated the
same way, except that each
time $\gamma$ hits a hexagon $f$ whose color has yet to be determined,
we perform a simple random walk
on the space of hexagons, beginning at $f$, and let $f$ assume the
color of the first black or
white hexagon hit by that walk. (See Figure~\ref{hesetuphexagons}.) In
other words, we color $f$
black with probability equal to the value at~$f$ of the function which
is equal to $1$ on the
black faces and $0$ on the white faces, and is {\it discrete harmonic}
at the undetermined faces
(i.e., its value at each such face is the mean of the values on the six
neighboring faces).

%fig1
%f1 ###

Denote by $h_n$ the value of this function after $n$ steps of the
harmonic explorer process;
that is, $h_n(f)$ is $1$ if $f$ is black, $0$ if $f$ is white, and
discrete harmonic on the faces
of undetermined color. Note that $h_n(f)$ is also the probability that
a random walk on faces,
started at $f$, hits a black face before hitting a white face. It is
easy to see (and
proved below) that for any fixed $f$, $h_n(f)$ is a martingale---and
that the harmonic explorer
is the only random path with this property. We will see later that
$\mathrm{SLE}_4$ is the only
random path with a certain continuous analog of this property.

% The stochastic Loewner evolution (SLE) processes were introduced in~
% for the scaling limits of statistical physics processes in two
%dimensions. The SLE is obtained by
% running Loewner's differential equation with Brownian motion as the
%driving parameter.
% Though we briefly recall the definitions below, the reader is advised
%to consult the
% survey~\cite{\WernerSurvey} for background on SLE.

It was conjectured in~\cite{16} and proved in  \cite{18} that if the
interior hexagons are each colored via critical Bernoulli percolation
(i.e., $p=p_c=1/2$), then, in a certain well-defined sense, the random
paths $\gamma$ tend to the stochastic Loewner evolution with parameter
$\kappa= 6$ ($\mathrm{SLE}_6$) as the hexagonal mesh gets finer. (See
the survey~\cite{19} for background on SLE.) It has been further
conjectured~\cite{13} that if colors are instead chosen from a critical
FK cluster model (where one weights configurations according to the
total number of clusters and the lengths of their interfaces), then
$\gamma$ will converge to some ${\mathrm{SLE}}_\kappa$ with $4 <
\kappa< 8$, where $\kappa$ depends on the weight parameters. We will
prove that, as the mesh gets finer, the harmonic explorer converges to
chordal $\mathrm{SLE}_4$.

There are also natural variants of the harmonic explorer; for example,
we might replace the
honeycomb lattice with another three-regular lattice or replace the
simple random walk on faces
with a different periodic Markov chain. One may even use a
non-three-regular lattice provided one
fixes an appropriate ordering (say, left to right) for determining the
color of multiple
undetermined faces that are ``hit'' simultaneously by the HE path.
Provided the simple random walk
converges to Brownian motion as the mesh gets finer, we see no barrier
to extending our results to
all of these settings. Our proofs are more like the LERW proofs
in~\cite{9} (which hold for
general lattices) than the percolation proof in \cite{18}
(which uses the invariance of
the lattice under $2\pi/3$ rotation in an essential way).
However, for simplicity, we will focus only on the hexagonal lattice in
this paper.

Although physicists and mathematicians have conjectured that many
models for random self-avoiding
lattice walks have conformally invariant scaling limits
[e.g., the
infinite self-avoiding walk,
critical percolation cluster boundaries on two-dimensional lattices,
critical Ising model interfaces, critical FK cluster boundaries and
$O(n)$ model strands, etc.],
rigorous proofs are available only in the following cases: percolation
interface on
the hexagonal lattice (which converges to
chordal~$\mathrm{SLE}_6$), harmonic
explorer (chordal $\mathrm{SLE}_4$),
loop erased random walk (LERW) on a periodic planar graph (radial
${\mathrm{SLE}}_2$ \cite{9}), the
uniform spanning tree (UST) boundary (chordal
${\mathrm{SLE}}_8$  \cite{9}) and the boundaries of
simple random walks (essentially ${\mathrm{SLE}}_{8/3}$
\cite{8}---here conformal invariance
follows easily from
the conformal invariance of Brownian motion).

The harmonic explorer is similar in spirit to the loop erased random
walk (LERW) and diffusion
limited aggregation (DLA). All three models are processes based on
simple random walks, and their
transition probabilities may all be computed using discrete harmonic
functions with appropriate
boundary conditions. Since simple random walks on
two-dimensional
lattices have a conformally
invariant scaling limit (Brownian motion), and since harmonicity (in
the continuous limit) is a
conformally invariant property, one might expect that all three models
would have conformally
invariant scaling limits. However, simulations suggest that DLA is not
conformally invariant.

This paper follows the strategy of~\cite{9}, and uses some of
the techniques from that
paper. We will freely quote results from~\cite{9}, and therefore
advise the reader to have
a copy of~\cite{9} on hand while reading the present paper.

%s2 ###
\section{A martingale property of chordal
$\mathrm{SLE}_4$.}

The purpose of this section is to briefly review some background about Loewner's
equation and SLE, and then present the basic strategy of the paper.
For more details, the reader is encouraged to
consult~\cite{19}
or~\cite{7}.

Let $T>0$.
Suppose that $\gamma\dvtx [0,T]\to\overline{\H}$
is a continuous simple path
in the closed
upper half plane $\overline{\H}$ which satisfies
$\gamma[0,T]\cap\R=\{\gamma(0)\}=\{0\}$.
For every $t\in[0,T]$, there is a unique conformal homeomorphism
$g_t\dvtx \H\setminus\gamma[0,t]$ which satisfies the so-called
{\it hydrodynamic}
normalization at infinity
\[
\lim_{z\to\infty} g_t(z)-z=0.
\]
The limit
\[
\capacity(\gamma[0,t]):=\lim_{z\to\infty} z\bigl(g_t(z)-z\bigr)/2
\]
is real and monotone increasing in $t$.
It is called the (half plane) {\it capacity} of $\gamma[0,t]$ from
$\infty$, or
just capacity, for short.
Since $\capacity(\gamma[0,t])$ is also continuous in $t$,
it is natural to reparameterize $\gamma$ so that
$\capacity(\gamma[0,t])=t$.
Loewner's theorem states that in this case the maps $g_t$ satisfy
his differential equation
%
%e2.1 ###
\begin{equation}\label{e.chordal}
\p_t g_t(z) =
\frac{2}{g_t(z)-W(t)},\qquad g_0(z)=z,
\end{equation}
where $W(t)=g_t(\gamma(t))$.
(Since $\gamma(t)$ is not in the domain of definition of $g_t$, the expression
$g_t(\gamma(t))$ should be interpreted as a limit of
$g_t(z)$ as $z\to\gamma(t)$ inside $\H\setminus\gamma[0,t]$. This limit
does exist.) The function $W(t)$ is
continuous in $t$, and is called the {\it driving parameter} for $\gamma$.
% In particular, Loewner's theorem tells us that the path $\gamma$ can
% be reconstructed from its driving parameter.

One may also try to reverse the above procedure.
Consider the Loewner evolution defined by the ODE~\eref{e.chordal},
where $W(t)$ is a continuous,
real-valued function.
The path of the evolution is defined
as $\gamma(t)=\lim_{z\to W(t)}g_t^{-1}(z)$,
where $z$ tends to $W(t)$ from within the upper half plane $\H$,
provided that the limit exists.
The process (chordal)
${\mathrm{SLE}}_\kappa$ in the upper half plane, beginning at $0$ and ending at
$\infty$, is the path
$\gamma(t)$ when $W(t)$ is $\sqrt\kappa B_t$, where $B_t=B(t)$ is a
standard one-dimensional Brownian motion.
(``Standard'' means $B(0)=0$ and $\E[B(t)^2]=t$, $t\ge0$.
Since $(\sqrt{\kappa} B_t \dvtx  t\ge 0)$ has the same distribution as
$(B_{\kappa t}\dvtx t\ge0)$,
taking $W(t)=B_{\kappa t}$ is equivalent.)
In this case a.s. $\gamma(t)$ does exist and is a continuous path.
See~\cite{13} ($\kappa\ne 8$)
and~\cite{9} ($\kappa=8$).

Fix $\kappa>0$, and assume now that $W(t)=\sqrt{\kappa} B_t$ and
$\gamma$ is ${\mathrm{SLE}}_{\kappa}$. Write $X = X(t,z) = g_t(z) -
W(t)$. Then $\arg X(t,z)$ gives $\pi$ times the probability that a
two-dimensional Brownian motion starting at $z$ first exits
$\H\setminus \gamma[0,t]$ either in $(-\infty,0)$ or on the left-hand
side of $\gamma[0,t]$. (This follows from conformal invariance of a
Brownian motion, run until its first exit point, and the fact that the
probability that a Brownian motion started at $z \in\mathbb H$ first
hits $\mathbb R$ at $(-\infty, 0]$ is $\arg(z)/\pi$. The latter fact
may be seen by conformally mapping the half plane to a strip using the
function $z \rightarrow\log(z)$.) In other words, for fixed $t$, $\arg
X(t,z)$ is the harmonic function that is equal to $\pi$ on one side of
$\gamma[0,t]$ and $0$ on the other. For short, we will sometimes refer
to the quantity $\arg X(t,z)$ as simply the {\it angle} of $z$ at time
$t$.

Now, using It\^{o}'s formula, we compute
the It\^{o} derivatives of $X$ and $\log X$:
\begin{eqnarray*}
d X &=&  \frac{2}{X}\,dt - \sqrt{\kappa}\,dB_t,
\\
d \log X &=&  \frac{2}{X^2}\,dt - \frac{\sqrt{\kappa}}{X}\,dB_t - \frac
{\kappa}{2X^2}\,dt = \frac{4-k}{2X^2}\,dt -
\frac{\sqrt{\kappa}}{X}\,dB_t.
\end{eqnarray*}

When $\kappa= 4$, we have $d\log X = -{2}{X^{-1}}\,dB_t$, and hence
$d \arg X =\break -
\Im({2}{X^{-1}})\,dB_t$. In particular, this implies that for any fixed
value of $z$,\break $\arg
X(t,z)/\pi$ is a martingale which is bounded in the interval $[0,1]$.
The value of this martingale
a.s. tends to  either zero or $1$ as $t$ tends to infinity, depending on
whether $z$ is on the left or the
right side of the path (see~\cite{15}, Lemma 3).
Hence, at a
fixed time $t$, $\arg X(t,z)/\pi$
represents the {\em probability} that, conditioned on the SLE path up
until time $t$, the point $z$
will lie to the left of the path.

It is easy to see (and shown below) that a discrete version of this
property holds for the harmonic
explorer. The strategy of our $\mathrm{SLE}_4$ proof will be, roughly speaking,
to show that the fact that
this property holds at two distinct values of $z$ is enough to force
the Loewner driving process
for the path traced by the harmonic explorer to converge to Brownian
motion. This is because the
fact that $\arg X(t,z)$ is a martingale at $z$ gives a linear
constraint on the drift and diffusion
terms at that point, and using two values of $z$ gives two linear
constraints, from which it is
possible to calculate the drift and diffusion exactly. The arguments
and error bounds needed to
make this reasoning precise are essentially the same as those
given in~\cite{9} (but the
martingales considered there are different). The fact that the Loewner
driving process converges to
Brownian motion will enable us to conclude that HE converges to
$\mathrm{SLE}_4$ in the Hausdorff topology.
We will then employ additional arguments to show that the convergence
holds in a stronger topology.

We remark that we will reuse this strategy in \cite{17} to prove
that a certain zero level
set of the discrete Gaussian free field (defined on the vertices of a
triangular lattice, with
boundary conditions equal to an appropriately chosen constant $\lambda$
on the left boundary and
$-\lambda$ on the right boundary) converges to chordal $\mathrm{SLE}_4$. To
keep notation consistent
with~\cite{17} (which will cite the present paper), we will use the dual
formulation (representing hexagons by vertices of the triangular
lattice) in our precise statements
and proofs below.

%s3 ###
\section{Statements of main results.} \label{s.HEsetup}

%s3.1 ###
\subsection{Notation and basic properties of HE.}\label{s.HEconstruction}

We now introduce the precise combinatorial notation for HE that we will
use in our proofs. First,
the triangular grid in the plane will be denoted by
$\mathit{TG}$. Its
vertices, denoted by $V(\mathit{TG})$,
are the sublattice of $\C$ spanned by $1$ and $e^{2\pi i/6}$; two
vertices are adjacent if
their difference is a sixth root of unity. If $D\subset\C$, and
$z\in\C$, let
$\inr z(D)$ denote the {\it inradius} of $D$ about $z$; that is,
$\inr z(D):=\inf\{|w-z|\dvtx w\notin
D\}$. Let $\Doms$ denote the set of domains $D\subset\C$ % with $0\in D$
whose boundary is a simple closed curve which is a union of edges from
the lattice $\mathit{TG}$.

If $V_0$ is any set of vertices in $V(\mathit{TG})$, and $h\dvtx V_0\to\R$
is a
bounded function, then there
exists a unique bounded function $\bar{h}\dvtx
V(\mathit{TG})\to\R$ which agrees
with $h$ in $V_0$ and is
harmonic at every vertex in $V(\mathit{TG})\setminus V_0$. This function is
called the {\it discrete
harmonic extension} of $h$.
[In fact, $\bar h(v)$ is the expected value
of $h$ at the point at
which a simple random walk started at $v$ hits $V_0$. Uniqueness is
easily established using the
maximum principle.]

\begin{figure}

\includegraphics{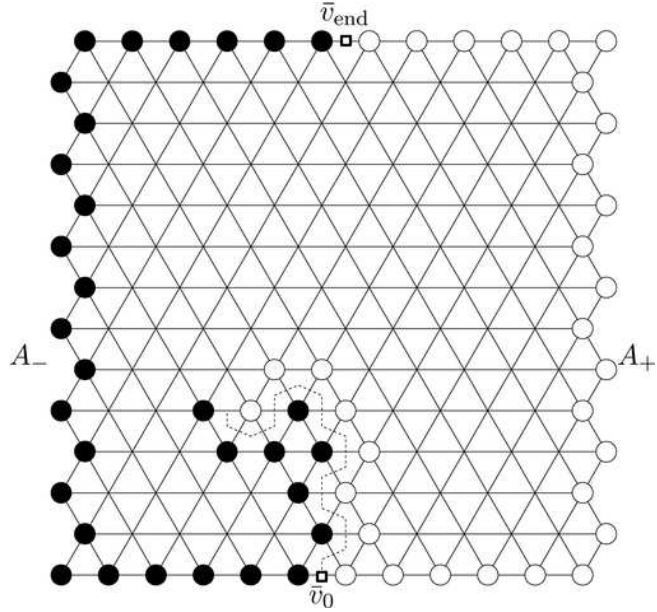}

\caption{A dual perspective on
Figure \up{\protect\ref{hesetuphexagons}}.}
\label{f.setup}
\end{figure}

Let $D\in\Doms$.
Let $V_0:=V(\mathit{TG})\cap\p D$ denote the set of vertices
in $\p D$. Let $\oo_0$~and~$\bo$ be the centers of two distinct edges of the grid $\mathit{TG}$ on $\p D$.
(See Figure~\ref{f.setup}.) Let $\AA_+$ (resp.
$\AA_-$) be the positively (resp. negatively) oriented arc of
$\p D$ from $\oo_0$ to $\bo$.
Define $\HH_0\dvtx V_0\to\{0,1\}$ to be $1$ on $V_0\cap\AA_+$, and $0$ on
$V_0\cap\AA_-$. The HE
[depending on the triple $(D,\oo_0,\bo)$] is a random simple path
from~$\oo_0$ to~$\bo$ in
$\closure D$. Let $X_1,X_2,\ldots$ be i.i.d. random variables, uniform
in the interval $[0,1]$.
(These will be the ``coin flips'' needed to generate the HE.) Let
$T_1\subset D$ be the triangle of
$\mathit{TG}$ whose boundary contains $\oo_0$
let $v_1$ be the vertex of $T_1$
that is not on the edge
containing $\oo_0$ and let $V_1:=V_0\cup\{v_1\}$. Let $v_1^{\prime}$ be the
middle of the edge of $T_1$ which is on the positively oriented arc
from $\oo_0$ to $v_1$, and let
$v_1^{\prime\prime}$ be the middle of the edge of $T_1$ which is on the positively
oriented arc from $v_1$ to
$\oo_0$. Let $p_1$ be the value at $v_1$ of the discrete harmonic
extension of $\HH_0$. If $X_1\le
p_1$, we let $\oo_1:=v_1^{\prime\prime}$, and otherwise
$\oo_1:=v_1^{\prime}$. The beginning
of the HE path is chosen as
the union of the two line segments from $\oo_0$ to the center of the
triangle $T_1$ and then to
$\oo_1$. Now define $\HH_1 \dvtx
V_1\to\{0,1\}$ to equal $\HH_0$ on $V_0$ and
set $\HH_1(v_1):=\mathbh{1}_{X_1\le
p_1}$ if $v_1\notin V_0$. This defines the first step of the HE.

%fig2
%f2 ###

The process continues inductively. Assuming that $n\ge1$ and
$\oo_n\notin\p D$, let $T_{n+1}$ be the triangle of $\mathit{TG}$
containing $\oo_n$ but not $\oo_{n-1}$. Let $v_{n+1}$ be
the vertex of $T_{n+1}$ which is not on the edge containing $\oo_n$,
and let $V_{n+1}:=V_n\cup\{v_{n+1}\}$.
Let $p_{n+1}$ be the value at $v_{n+1}$ of the discrete harmonic extension
of $\HH_n$.
Let $v_{n+1}^{\prime}$ and $v_{n+1}^{\prime\prime}$ be the two midpoints of edges of
$T_{n+1}$ that lie on the positively oriented
arcs of $\p T_{n+1}$ from $\oo_n$
to $v_{n+1}$ and from $v_{n+1}$ to $\oo_n$, respectively.
If $X_{n+1}\le p_{n+1}$ let $\oo_{n+1}:= v_{n+1}^{\prime\prime}$,
and otherwise $\oo_{n+1}:=v_{n+1}^{\prime}$. Let the next step of the HE
consist of the segments from $\oo_n$ to the center of~$T_{n+1}$
and from the center of~$T_{n+1}$ to~$\oo_{n+1}$.
Also, let $\HH_{n+1}$ agree with $\HH_n$, where $\HH_n$ is defined
and set $\HH_{n+1}(v_{n+1}):= \mathbh{1}_{X_{n+1}\le p_{n+1}}$ if
$v_{n+1}\notin V_n$.

It is easy to verify that this procedure a.s. terminates
when $\oo_n=\bo$, and that the HE so defined is a simple
path from $\oo_0$ to $\bo$. Let $N$ denote the termination
time; that is, the $n$ such that $\oo_n=\bo$.

\begin{lemma}\label{l.hmart}
Let $h_n$ denote the discrete harmonic extension of $\HH_n$, and
let $v\in V(\mathit{TG})\cap D$. Then $h_n(v)$ is a martingale
and $h_N(v)\in\{0,1\}$.
\end{lemma}

\begin{pf}
Given $X_1,\ldots,X_n$, we have
$h_{n+1}(v_{n+1})=1$ with probability\break
$h_n(v_{n+1})$ and otherwise $h_{n+1}(v_{n+1})=0$.
Consequently, $\E[h_{n+1}(v_{n+1})| X_1,\ldots,\break X_n]=h_n(v_{n+1})$.
Note that $h_{n}$ is also the discrete harmonic extension of
its restriction to $V_{n+1}$, and similarly for $h_{n+1}$.
Since the harmonic extension is a linear operation and
$\E[h_{n+1}(v)| X_1,\ldots,X_n]=h_n(v)$ for $v\in V_{n+1}$,
the same relation holds for every $v$. Thus $h_n(v)$ is
a martingale. The claim that $h_N(v)\in\{0,1\}$ is clear.
\end{pf}

\begin{rem}\label{r.HEmark}
The evolution of the HE path may be viewed as a Markov chain on the collection
of appropriately marked domains.
At the $n$th step, the chain is at $(D_n,\oo_n,\bo)$,
where $D_n$ is the connected component of
$D\setminus\bigcup_{j=1}^n T_j$
that has $\bo$ on its boundary.
\end{rem}

%s3.2 ###
\subsection{Convergence of HE\up{:} statement.}

Let $D\in\Doms$. We assume the setup and notation of Section
\ref{s.HEconstruction}. Let $\gamma\dvtx [0,N]\to D\cup\{\oo_0,\bo\}$
be the HE path with the parameterization proportional to arclength,
where $\gamma(n)=\oo_n$ for $n\in\{0,1,\ldots,N\}$. Let $\phi \dvtx
D\to\H$ be a conformal map onto $\H$ that takes $\oo_0$ to $0$ and
$\bo$ to $\infty$. Note that $\phi$ is unique up to positive scaling,
$\phi(\AA_+)=(0,\infty)$ and $\phi(\AA_-)=(-\infty,0)$. Let
$\bb:=\phi^{-1}(i)$.

Instead of rescaling the grid, we consider larger and larger domains
$D$. The quantity
$\rr=\rr(D,\phi):=\inr\bb(D)$ turns out to be the appropriate
indicator of the size of $D$, from
the perspective of the map $\phi$. Indeed, if $\rr$ is small, then the
image under $\phi$ of the
grid $\mathit{TG}$ in $D$ is not fine near $i$, and we cannot expect $\phi\circ
\gamma$ to look like
$\mathrm{SLE}_4$. As we will see,
$\phi\circ\gamma$ does approach $\mathrm{SLE}_4$ when
$\rr\to\infty$. Let $\hg$ be
the path $\phi\circ\gamma$, parameterized by capacity from $\infty$ in
$\closure\H$, and let $\tg$
be the $\mathrm{SLE}_4$ path in~$\closure\H$.

Let $d_*(\cdot,\cdot)$ be the metric on $\closure\H\cup\{\infty\}$
given by $d_*(z,w)=|\Psi(z)-\Psi(w)|$, where $\Psi(z):=(z-i)/(z+i)$
maps $\closure\H\cup\{\infty\}$ onto $\closure\U$.
If $z\in\closure\H$, then $d_*(z_n,z)\to 0$ is equivalent
to $|z_n-z|\to0$, and $d_*(z_n,\infty)\to 0$ is equivalent
to $|z_n|\to\infty$.

Note that although we started by mapping our domain
$D$ to the half plane (with boundary points $0$, $\infty$ and inradius
measured from the
preimage of $i$), the above metric corresponds to a mapping to the unit
disc (with
boundary points $-1$, $1$ and inradius measure from the preimage of
$0$). The half plane
is the most convenient setting for describing Loewner evolution and
chordal $\mathrm{SLE}_4$, but
the metric derived from the unit disc map is more convenient because it
is compact.

\begin{thm}\label{t.heunifconv}
As $\rr\to\infty$, the law of $\hg$ tends to the
law of the $\mathrm{SLE}_4$ path~$\tg$, with respect to uniform convergence in
the metric
$d_*$. In other words, for every $\eps>0$ there is some $R=R(\eps)$
such that if $\rr>R$, then there is a coupling of
$\hg$~and~$\tg$ such that
\[
\P\bigl[\sup\bigl\{d_*\bigl(\tg(t),\hg(t)\bigr)\dvtx
t\in(0,\infty)\bigr\}>\eps\bigr]<\eps.
\]
\end{thm}

%s4 ###
\section{The driving process converges to BM.}

Let $W=W(t)$ denote
the Loewner driving process for $\hg$. Let $B\dvtx [0,\infty)\to\R$ be a standard
one-\break dimensional Brownian motion.
A slightly weaker form of Theorem~\ref{t.heunifconv} will follow as a
consequence of the fact
that for every $T>0$ the restriction of $W$ to $[0,T]$ converges in law
to the restriction of
$t\mapsto 2 B(t)$ to $[0,T]$. This, in turn, will be a consequence of
the following local
statement.

\begin{proposition}\label{p.helocal}
For $n\in[0,N]$ let $t_n:=\capacity(\psi\circ\gamma[0,n])$,
$\DD_n:=D\setminus\gamma[0,n]$, and let $\phi_n\dvtx \DD_n\to\H$ be the
conformal map normalized by $\phi_n\circ\phi^{-1}(z)-z\to0$ as
$z\to\infty$ in $\H$. For every $\delta\in(0,1)$ there is an
$R=R(\delta)>0$ such that the following holds. Fix any $n\in\N$. On the
event $\ev A_1=\ev A_1(n):=\{n<N\}$, let $m$ be the least integer
larger than $n$ such that
$\max\{t_m-t_n,(W(t_m)-W(t_n))^2\}\ge\delta^2$. \up{(}Note that $m\le
N\!$, since $t_N=\infty$.\up{)} Set $\aa:= \phi_n^{-1}(i+W(t_n))$ and
let $\ev A_2=\ev A_2(n)$ be the event $\{\inr\aa(D)\ge R\}$. Then
%
%e4.1 ###
\begin{equation}\label{e.hemart}
\E \bigl[ W(t_m)\md\gamma[0,n]\bigr]
= W(t_n)+O(\delta^3)
\end{equation}
and
%
%e4.2 ###
\begin{equation}\label{e.hevar}
\E \bigl[\bigl(W(t_m)-W(t_n)\bigr)^2\md\gamma[0,n]\bigr] =
4\E \bigl[t_m-t_n\md\gamma[0,n]\bigr] +O(\delta^3)
\end{equation}
both hold on the event $\ev A_1\cap\ev A_2$.
\end{proposition}

Here, and below, $O(f)$ represents any quantity whose
absolute value is bounded by $cf$, where $c$ is any fixed constant.

The strategy for proving the proposition is as follows.
We use Lemma~\ref{l.hmart} to conclude that
$\E[h_m(v)-h_n(v)\md\gamma[0,n]]=0$.
Since $h_j(v)$ is discrete-harmonic, it is approximately equal to
the harmonic function on $\DD_j$ with the corresponding boundary values.
The difference $h_m(v)-h_n(v)$ can then be approximated by a function
of $t_m-t_n$ and $W(t_m)-W(t_n)$. Applying
$\E[h_m(v)-h_n(v)\md\gamma[0,n]]=0$
for two distinct choices of $v$ then
gives the relations~\eref{e.hemart} and~\eref{e.hevar}.

We start with a lemma describing the
approximation of $h_j$ by a (nondiscrete)
harmonic function.

\begin{lemma} \label{l.harmconv}
Given any $\eps>0$ there is an $r=r(\eps)>0$ such that
for every vertex $v\in V(\mathit{TG})$ and every $j<N$,
if $\inr v(\DD_j)>r$, then
%
%e4.3 ###
\begin{equation}\label{e.harmconv}
\bigl|h_j(v)- \tilde{h}\bigl(\phi_j(v)-W(t_j)\bigr)\bigr|<\eps,
\end{equation}
where $\tilde h(z)=1-(1/\pi)\arg z$.
\end{lemma}

\begin{pf}
Note that $\tilde h \dvtx \H\to(0,1)$ is harmonic
and has the boundary values
$0$ on $(-\infty,0)$ and $1$ on $(0,\infty)$.
Since $W(t_j)=\phi_j(\gamma(j))$, $z\mapsto\tilde h(\phi
_j(z)-W(t_j))$
is harmonic in $\DD_j$, and has boundary values $0$ on $\AA_-$
and on the
``left side'' of $\gamma[0,j]$ and
$1$ on $\AA_+$ and the ``right side'' of $\gamma[0,j]$.
Since $h_j$ is a discrete harmonic function with similar boundary conditions,
the statement of the lemma can be obtained as a consequence of the convergence
of random walk on $\mathit{TG}$ to Brownian motion.
We leave the details to the reader.
(Also note that more delicate but similar estimates are
given in~\cite{9}, Section~\ref{sec5}.)
\end{pf}

\begin{pf*}{Proof of Proposition~\ref{p.helocal}}
Assume $\ev A_1$. We claim that there is an absolute constant
$\delta_0>0$ such that $\inr\aa(\DD_m)\ge\frac12\inr\aa(\DD_n) - 1$ if
$\delta<\delta_0$. Let $z$ be on the circle $|z|=\inr\aa(\DD_n)/2$.
Since $\Im\phi_n(\aa)=1$, the Koebe distortion theorem implies a
positive constant lower bound for $\Im\phi_n(z)$ (see, e.g., \cite{12},
Section~1.3). Let $g_t$ be the Loewner chain driven by $W(t)$. Then
$\phi _j=g_{t_j}\circ\phi$. By Loewner's equation~\eref{e.chordal},
$\frac{d}{dt}\Im g_t(z)\ge-2 /\Im g_t(z)$, which implies
$\frac{d}{dt}(\Im g_t(z))^2\ge-4$. Thus, $\tau(z)\ge t_n+(\Im\phi
_n(z))^2/4$. Since $t_{m-1}-t_n\le\delta^2$, it follows that
$z\notin\gamma[0,m-1]$ if $\delta<\delta_0$, where $\delta_0$ is the
infimum of all possible values for $\Im\phi_n(z)/2$. In that case,
$\inr\aa(\DD_{m-1})\ge\inr\aa(\DD_n)/2$, which implies our claim $\inr
\aa(\DD_m)\ge\frac 12\inr\aa(\DD_n) - 1$. We will henceforth assume,
with no loss of generality, that $\delta<\delta_0$. Note that the above
argument also gives a positive lower bound on
$\Im\phi_{m-1}(\aa)$.

Now fix some vertex $w_0\in\DD_n\cap V(\mathit{TG})$ satisfying
$|w_0-\aa|< \inr\aa(\DD_n)/6$. Let $R$ be larger than $100\max\{
1,r(\delta^3)\}$,
in the notation of Lemma~\ref{l.harmconv}.
Assume now that $\ev A_2$ holds.
Then we may apply~\eref{e.harmconv} with $j=n,m$, $v=w_0$ and $\eps
=\delta^3$.
Since $h_j$ is a martingale, it satisfies
$\E[h_m(w_0)\md\gamma[0,n]]
= h_n(w_0)$,
and so we get from~\eref{e.harmconv}
%
%e4.4 ###
\begin{equation}\label{e.tr}
\E\bigl[\tilde h\bigl(\phi_m(w_0)-W(t_m)\bigr)\md\gamma[0,n]\bigr]
= \tilde{h}\bigl(\phi_n(w_0)-W(t_n)\bigr)
+O(\delta^3).
\end{equation}

Below, we need the relations
%
%e4.5 ###
\begin{equation}\label{e.aprio}
\forall\, t\in[t_n,t_m] \qquad |W(t)-W(t_n)|= O(\delta),\qquad
t_m-t_n=O(\delta^2).
\end{equation}
By our choice of $m$, we have the first relation when
$t=t_n,t_{n+1},\ldots,t_{m-1}$
and the second relation when $t_{m-1}$ replaces $t_m$.
The relations~\eref{e.aprio} will follow by assuming that $R$ is large enough.
Indeed, if $j\in\{n,\ldots,m-1\}$ and $R$ is large, then the harmonic
measure from $\aa$ of $\gamma[j,j+1]$ in $\DD_{j}$ is
$O(\delta)$. (The Beurling projection theorem \cite{1} tells us
that $R>\delta^{-2}$
suffices.) By conformal invariance of harmonic measure,
the harmonic measure from $\phi_j(\aa)$ of $\phi_j\circ\gamma[j,j+1]$
in $\H$ is $O(\delta)$.
We want to use this to conclude that $\diam(\phi_j\circ\gamma
[j,j+1])=O(\delta)$.
Note that $\phi_j(\aa)=g_{t_j}\circ\phi(\aa)$. Above, we have seen that
there is a constant positive lower bound for $\Im\phi_{m-1}(\aa)$.
By~\eref{e.chordal}, $\Im g_t(z)$ is monotone
decreasing in $t$. Hence, $\Im g_t\circ\phi(\aa)$ has a constant
positive lower bound
for $t\le t_{m-1}$. By~\eref{e.chordal}, we get $ |\p_t( g_t\circ
\phi(\aa))|=O(1)$
for $t\le t_{m-1}$. Integrating then gives $|\phi_j(\aa)-\phi_n(\aa
)|\le O(\delta^2)$
for $j=n,n+1,\ldots,m-1$. As $W_j=\phi_j( \oo_j)\in\phi_j\circ\gamma[j,j+1]$
the distance from $\phi_j(\aa)$ to $\phi_j\circ\gamma[j,j+1]$ is
$O(1)$. Consequently,
the harmonic measure estimate gives the bound $\diam\phi_j\circ\gamma
[j,j+1]=O(\delta)$.
The needed estimates~\eref{e.aprio} now follow
from~\cite{9}, Lemma 2.1, since
$\phi_j(\gamma[j,j+1])$ is the set of points hitting the real line
under Loewner's
evolution~\eref{e.chordal} in the time interval $[t_j,t_{j+1}]$.

Let $z_t:=g_t\circ\phi(w_0)$.
Since we have $\phi_j(w_0)=z_{t_j}$, we may obtain $\phi_m(w_0)=z_{t_m}$
from $\phi_n(w_0)=z_{t_n}$ by flowing according
to Loewner's equation~\eref{e.chordal} between the times $t_n$ and $t_m$.
As before, we get the bound $|z_t-z_{t_n}|=O(\delta^2)$ for $t\in[t_n,t_m]$.
Since $|W(t)-W({t_n})|=O(\delta)$,
we have
\[
\frac2{z_t-W(t)}=\frac2{z_{t_n}-W(t_n)}+ O(\delta),\qquad t\in
[t_n,t_m].
\]
By integrating this relation over $[t_n,t_m]$,~\eref{e.chordal} gives
%
%e4.6 ###
\begin{equation}\label{e.dz}
z_{t_m}-z_{t_n}=
\phi_m(w_0)-\phi_n(w_0)= \frac{2(t_m-t_n)}{\phi_n(w_0)-W(t_n)} +
O(\delta^3).
\end{equation}

Consider now $F(z,W):=\tilde h(z-W)$. We want an estimate for
\[
F\bigl(z_{t_m},W(t_m)\bigr)
=\tilde{h} \bigl(\phi_m(w_0)-W(t_m)\bigr)
\]
up to $O(\delta^3)$ terms.
For that purpose, we use a Taylor expansion of $F$ about $(z_{t_n},W(t_n))$.
Since $z_{t_m}-z_t=O(\delta^2)$ and $W(t_m)-W(t_n)=O(\delta)$,
it suffices to take the terms up to the first derivative of $F$ with respect
to $z$ and the second derivative of $F$ with respect to $W$, and no
mixed terms.
Hence,
\begin{eqnarray*}
&&
\tilde{h}
\bigl(\phi_m(w_0)-W(t_m)\bigr)
-\tilde{h}\bigl(\phi_n(w_0)-W(t_n)\bigr)
\\
&&\qquad  =
\p_z F_{(z_{t_n},W(t_n))} \bigl(z_{t_m}-z_{t_n}\bigr)+
\p_W F_{(z_{t_n},W(t_n))} \bigl(W(t_m)-W(t_n)\bigr)\\
&&\qquad \quad{}+ \tfrac{1}{2} \p_W^2 F_{(z_{t_n},W(t_n))}
\bigl(W(t_m)-W(t_n)\bigr)^2 + O(\delta^3).
\end{eqnarray*}
(Since $z$ is complex, $\p_z F_{(z_{t_n},W(t_n))}$ is actually
a linear map from $\C$ to $\R$.)
By~\eref{e.tr}, the conditional expectation of the left-hand side given
$\gamma[0,n]$ is $O(\delta^3)$. After calculating the derivatives and
applying~\eref{e.dz}, we get
%
%e4.7 ###
\begin{equation}
\label{e.gettingthere}
\quad\begin{array}{rl}
\ds O(\delta^3) =  & \ds
2 \Im\bigr(\bigl(\phi_n(w_0)-W(t_n)
\bigr)^{-2}\bigr) \E \bigl[t_m-t_n\md\gamma[0,n]\bigr]
\\
\noa
& \ds {}- \Im\bigl(\bigl(\phi_n(w_0)-W(t_n)\bigr)^{-1}\bigr)
\E\bigl[W(t_m)-W(t_n) \md\gamma
[0,n]\bigr]
\\
\noa & \ds {}- \tfrac12
\Im\bigl(\bigl(\phi_n(w_0)-W(t_n)\bigr)^{-2}\bigr) \E\bigl[
\bigl(W(t_m)-W(t_n)\bigr)^2\md\gamma [0,n]\bigr].
\end{array}
\end{equation}
We now assume that $R>\delta^{-2}$. The Koebe distortion theorem
(again, see \cite{2}, Section~1.3) then implies
that a vertex $w_1\in V(\mathit{TG})$ closest to $\aa$ satisfies
$|\phi_n(w_1)-i-W(t_n)|=|\phi_n(w_1)-\phi_n(\aa)|=O(\delta^2)$.
The Koebe distortion theorem also shows that a vertex $w_2\in V(\mathit{TG})$ closest
to $\phi_n^{-1}(i+W(t_n)+1/100)$ satisfies
$|\phi_n(w_2)-i-W(t_n)-1/100|=O(\delta^2)$
and $|w_2-\aa|<\inr\aa(D)/6$.
Consequently, we may apply~\eref{e.gettingthere}
with $w_0$ replaced by each of $w_1,w_2$.
With $w_0=w_1$, we get~\eref{e.hemart}. Now eliminating the term
$\E[W(t_m)-W(t_n) \md\gamma[0,n]]$
from~\eref{e.gettingthere}
[since it is $O(\delta^3)$]
and applying~\eref{e.gettingthere}
with $w_0=w_1$ gives~\eref{e.hevar}.
\end{pf*}

\begin{corollary}\label{c.drive}
Fix $T\ge1$. As $\rr\to\infty$,
the restriction of $t\mapsto W(t/4)$ to $[0,T]$ converges in law to the
corresponding
restriction of standard Brownian motion.
\end{corollary}

\begin{pf}
Let $\eps\in(0,1)$, and
let $\tilde T:=\sup\{t\in[0,T]\dvtx |W(t/4)|\le\eps^{-1}\}$.
Let $I:=\{n\in\N\dvtx  t_n\le\tilde T\}$.
In order to apply Proposition~\ref{p.helocal} at every $n\in I$, we
need to verify
$\ev A_1(n)\cap\ev A_2(n)$
for such $n$. From~\cite{9}, Lemma 2.1, we get that
$t_N=\infty$ or $\{W(t)\dvtx t\in[0,t_N]\}$ is unbounded, which implies
$\ev A_1(n)$ for $n\in I$.
Since $\phi_n(\aa)=i+W(t_n)$ and $g_{t_n}=\phi_n\circ\phi^{-1}$,
we have $g_{t_n}\circ\phi(\aa)=i+W(t_n)$.
We claim that there is a compact subset $K\subset\H$, which depends only
on $\eps$~and~$T$, such that $\phi(\aa)\in K$ holds for each $n\in I$. Indeed,
$g_t\circ\phi(\aa)$ flows according to~\eref{e.chordal} starting from
$\phi(\aa)$ at $t=0$ to $i+W(t_n)$ at $t=t_n$. For every $t\in
[0,t_n]$, we have
\mbox{$\Im g_t\circ\phi(\aa)\ge1$}, by the monotonicity of $\Im g_t$ with
respect to $t$.
By~\eref{e.chordal}, this shows that $|\p_t g_t\circ\phi(\aa)|=O(1)$.
Hence, the bound $|\phi(\aa)|\le1+|W(t_n)|+O(T)\le1+\eps^{-1}+O(T)$.
We may therefore take $K=\{z\in\C\dvtx \Im(z)\ge1,|z|\le O(T+\eps
^{-1})\}$.
Since $\phi(\aa)$ lies in a compact subset of $\H$, the Koebe
distortion theorem
implies that $\rr\le O(1)\,\inr\aa(D)$. Thus we may assume that $\ev A_2(n)$
holds for every $n\in I$ provided we take $\rr\ge R^{\prime}$ for some constant
$R'=R'(\eps,T,\delta)$.

Now the proof that the restriction of $t\mapsto W(t/4)$ to $[0,\tilde T]$
converges in law to the corresponding stopped Brownian motion follows from
the proposition and the Skorokhod embedding theorem,
as in~\cite{9}, Section~3.3.
Standard Brownian motion is unlikely to hit $\{-\eps^{-1},\eps^{-1}\}$
before time $T$ if $\eps$ is small. Thus, we obtain the corollary by taking
a limit as $\eps\searrow 0$.
\end{pf}

%s5 ###
\section{Local Hausdorff convergence to
$\mathrm{SLE}_4$.}\label{sec5}

Let $\dhaus(\cdot,\cdot)$ denote the Hausdorff distance;
that is, for two nonempty sets $A,B\subset\C$,
\[
\dhaus(A,B):= \max\biggl\{
\sup_{a\in A}\inf_{b\in B} |a-b|,
\sup_{b\in B}\inf_{a\in A} |b-a|\biggr\}.
\]

\begin{lemma}\label{l.parthaus}
For every $T\ge1$
and $\eps>0$ there is some $R=R(\eps,T)$ so that
if $\rr>R$, then there is a coupling of $\hg$
and $\tg$ so that
\[
\P \bigl[ \sup \{ \dhaus(\hg[0,t],\tg[0,t])\dvtx
{0\le t\le T}\}>\eps\bigr]<\eps.
\]
\end{lemma}

\begin{pf}
We know that $\tg$ is a simple path, from~\cite{13}.
Let $\tilde g_t$ be the $\mathrm{SLE}_4$ Loewner chain corresponding to $\tg$,
and let $B$ be the Brownian motion so that the driving process for
$\tilde g_t$ is $B(4t)$. Then $\tilde g_t$ is
obtained by solving~\eref{e.chordal} with
$W(t)$ replaced by $B(4t)$.

Let $S(T,\eps)$ denote the set of points in $\overline\H$
whose distance from $\tg[0,T]$ is~$\varepsilon$.
Let $s\in S(T,\eps)$. By continuity of solutions of differential equations,
there is some $\delta=\delta(s,B)>0$
such that if $W\dvtx [0,T]\to\R$ is
measurable and satisfies $\sup\{|W(t)-B(4  t)|\dvtx t\in[0,T]\}\le
\delta$,
then the Loewner chain corresponding to $W$ satisfies $\tau(s)>T$
(this is also easy to verify directly). Moreover, the same $\delta$
would apply to every $s'\in\closure\H$ sufficiently close to $s$.
By compactness of $S(T,\eps)$, there is some $\delta=\delta(B)>0$, which
would work for every $s\in S(T,\eps)$. Now, this $\delta$ is random, as it
depends on $B$, but it is a.s. positive.
Therefore, there is a nonrandom $\delta_0>0$,
depending only on $T$ and $\varepsilon$,
such that $\delta_0$ would work for $B$
with probability at least $1-\varepsilon/2$.

By Corollary~\ref{c.drive} (and the well-known relation between
convergence in law and  a.s. convergence
\cite{3}, Theorem 11.7.2),
when $\rr$ is sufficiently large,
we may couple~$\gamma$ with standard Brownian motion $B(t)$ so that
\[
\P\bigl[\sup \{ |W(t)-B(4t)|\dvtx
t\in[0,T]\}\ge\delta_0\bigr]<\eps/2,
\]
where now $W(t)$ is the driving process for $\hg$.
Consequently, when $\rr$ is sufficiently large,
\[
\P\bigl[ \hg[0,T]\cap S(T,\eps)=\varnothing\bigr] \ge 1-\eps.
\]
Since $\hg[0,T]$ is connected and contains $0$, when
$\hg[0,T]$ is disjoint from $S(T,\eps)$ every point in $\hg[0,T]$ is
within distance $\eps$ from $\tg[0,T]$.

Now consider a sequence of pairs $(D,\phi)$ such that $\rr\to\infty$.
For each such pair we take a coupling of the corresponding $\hg[0,T]$
and $B$ such that $ \sup\{ |W(t)-B(4t)|\dvtx t\in[0,T]\}\to0 $ in
probability. Fix some $t\in[0,T]$. Since the collection of probability
measures on the (compact) Hausdorff space of closed nonempty subsets of
$\closure\H\cup\{\infty\}$ is compact under convergence in law, by
passing to a subsequence, if necessary, we get a coupling of $B$ and a
Hausdorff limit $\Gamma_t$ of $\hg[0,t]$. By the above, $\Gamma_t$ is
contained in $\tg[0,T]$. Moreover, it is clearly connected. Note that
the Carath\'eodory kernel theorem~(\cite{12}, Theorem~1.8) implies that
the maps $g_t\dvtx \H\setminus\hg[0,t]\to\H$ converge to the normalized
conformal map from $\H\setminus\Gamma_t$ to~$\H$. Consequently, the
capacity of $\Gamma_t$ is $t$ a.s. Thus, we conclude that
$\Gamma_t=\tg[0,t]$ a.s. Since the limit does not depend on the
subsequence, it follows that $\tg[0,t]$ is~a.s.~the Hausdorff limit
$\Gamma_t$ of $\hg[0,t]$. As $t\in[0,T]$ is arbitrary, we conclude that
a.s. $\Gamma_t=\tg[0,t]$ for every rational $t$ in $[0,T]$. The lemma
now follows, since $\tg$ is continuous and
$\Gamma_t\supset\Gamma_{t^{\prime}}$ when $t>t^{\prime}\ge0$.
\end{pf}

\section{Improving the topology.}

Lemma~\ref{l.parthaus} gives some form of convergence of~$\hg$ to~$\tg$.
Our goal now is to improve the quality of the convergence, in two ways.
First, we want to show that the convergence is locally uniform
(i.e., uniform on compact intervals $[0,T]\subset[0,\infty)$).
Later, it will be shown that the convergence is uniform when we use the
metric $d_*$ on $\closure\H\cup\{\infty\}$.

To understand the issues here, we describe two examples where one form
of convergence holds and another fails. We start with an example
similar to one appearing in~\cite{9}, Section 3.4. Let $\eps>0$, $a_j:=
i \eps(1-j^{-1})$ and $b_j:= ij\eps$. Let $\alpha _\eps$ be the
polygonal path determined by the points $a_1,  b_1+\eps, a_2, b_2-\eps,
a_3, b_3+\eps,\break a_4, b_4-\eps$, etc. Let $\alpha_0$ be the path
$t\mapsto i t$, $t\ge0$, reparameterized by capacity. Then the path
$\alpha_\eps$ reparameterized by capacity converges to $\alpha_0$ in
the sense of Lemma~\ref{l.parthaus}.\vadjust{\goodbreak} Moreover, the
Loewner driving process for $\alpha_\eps$ converges locally uniformly
to the constant $0$, which is the driving process for $\alpha _0$.
However, one cannot reparameterize $\alpha_\eps$ so that
$\alpha_\eps\to\alpha_0$ locally uniformly.

To illustrate the second issue, consider the polygonal path
$\beta_\eps$ determined by the points $0,i \eps^{-1},i+\eps,\infty$,
where the last segment can be chosen as any ray from $i+\eps$ to $\infty$
in $\H$. Then $\beta_\eps$, reparameterized by capacity,
does converge
locally uniformly to $\alpha_0$. However, it does not converge uniformly
with respect to the metric $d_*$.

%s6.1 ###
\subsection{Discrete excursions.}

The purpose of this subsection is to develop a tool which will be handy
for proving some upper bounds on probabilities of rare events for the HE,
the discrete excursion measure.
It is a discrete analogue of the (two-dimensional)
Brownian excursion as introduced
in~\cite{10}, Section~2.4.
A slightly different variant of the
continuous Brownian excursion was studied in~\cite{8}.

Let $D$ be a domain in the plane whose boundary is a subgraph of the
triangular lattice $\mathit{TG}$. (We
work here with the simple random walk on $\mathit{TG}$, but the results apply
more generally to other walks
on other lattices.) Let $V_\p$ denote the set of vertices in $V(\mathit{TG})\cap
\p D$. A directed edge of
$\mathit{TG}$ is just an edge of $\mathit{TG}$ with a particular choice of orientation
(i.e., a choice of the
initial vertex). If $e=[u,v]$ is a directed edge
of $\mathit{TG}$, then $\rev (e)=[v,u]$ will denote the
same edge with the reversed orientation. Let $\inE=\inE(D)$ denote the
set of directed edges of
$\mathit{TG}$ whose interiors intersect $D$ and whose
initial vertex is in $V_\p
$. Let $\outE=\outE(D)$
denote the set of directed edges of $\mathit{TG}$ whose interiors intersect $D$
and whose terminal vertex
is in $V_\p$; that is, $\outE=\rev(\inE)$. Let $E_1\subset\inE$ and
$E_2\subset\outE$. For every
$v\in V_\p$, let $X^v$ be a simple random walk on $\mathit{TG}$ that starts at
$v$ and is stopped at the
first time $t\ge1$ such that $X^v(t)\notin D$. Let $\DX^v$ denote the
restriction of the law of
$X^v$ to those walks that use an edge of $E_1$ as the first step and
use an edge of $E_2$ as the
last step. (This is zero if $v$ is not adjacent to an edge in $E_1$,
and generally it is not a
probability measure.) Finally, let $\DX=\DX_{(D,E_1,E_2)}:=\sum_{v\in
V_\p} \DX^v$. This is a
measure on paths starting with an edge in $E_1$, ending with an edge in
$E_2$ and staying in $D$
in between. It will be called the discrete excursion measure from $E_1$
to $E_2$ in $D$. When
$E_2=\outE$, we will often abbreviate
$\DX_{(D,E_1)}=\DX_{(D,E_1,\scriptsize{\outE})}$.

\begin{lemma}\label{l.dxgreen}
Let $D$ be as above, and let $E_1\subset\inE$.
Fix $v\in V(\mathit{TG})\cap D$, and for every path $\omega$ let $n_v(\omega)$ be
the number of times $\omega$ visits $v$. Then
%
%e6.1 ###
\begin{equation}\label{e.dxgreen}
\int n_v(\omega)\, d\DX_{(D,E_1)}(\omega) =
H\bigl(v,\rev(E_1)\bigr),
\end{equation}
where $H(v,E)=H_D(v,E)$ denotes the probability that a simple
random walk started from $v$ will first exit $D$ through an edge
in $E$. In particular, \mbox{$\int n_v\,
d\DX_{(D,\scriptsize{\inE})}=1$}.
\end{lemma}

\begin{pf}
Let $(\Omega,\mu)$ denote the probability space of random walks
starting at $v$ and stopped
when they first exit $D$. For a pair $(\omega_1,\omega_2)\in\Omega^2$,
let $f(\omega_1,\omega_2)$
denote the reversal of $\omega_1$ followed by $\omega_2$. Then $f$ is a
map from $\Omega^2$ to the
support of $\DX=\DX_{(D,\scriptsize{\inE})}$. Clearly,
$\mu\times\mu(\{\omega_1,\omega_2\})=\DX(\{f(\omega_1,\omega_2)\})$. If
$\omega'$ is in $\Omega'$,
the support of $\DX$, then the cardinality of the preimage
$f^{-1}(\omega')\subset\Omega^2$ is
precisely $n_v(\omega')$. Consequently, we have
\[
1=\mu\times\mu(\Omega^2)=
\sum_{\omega'\in\Omega'}
|f^{-1}(\omega')| \DX(\{\omega'\}) =\int
n_v\,d\DX.
\]
This proves
the claim in the case $E_1=\inE$. The general case is similarly
established.
\end{pf}

\begin{corollary}\label{c.hitball}
Let $D$ and $E_1$ be as above, and let $v\in V(\mathit{TG})\cap D$.
Assume that $\p D$ is connected.
Let $B$ be the ball
centered at $v$ whose radius is $\frac{1}{2} \inr v(D)$, and let
$\Gamma_B$ be the set of paths that visit $B$. Then
\[
c^{-1}
H_D\bigl(v,\rev(E_1)\bigr)
<\DX_{(D,E_1)}(\Gamma_B)<c
H_D\bigl(v,\rev(E_1)\bigr)
\]
for some absolute constant $c$.
\end{corollary}

\begin{pf}
It is well known that there is an absolute constant $c$ with
$c^{-1}<G_D(w,v)<c$ for every $w\in B$ such that $|v-w|$ is
at least half the radius of~$B$, where $G_D$ is Green's function.
See, for example, \cite{9}, (3.5), where the radius of the
ball $B$ is different, but the same proof applies.
Consequently, given that a random walk hits $B$, its expected
number of visits to $v$ before exiting $D$ is between
$c^{-1}$ and $c$. Thus, the corollary follows from~\eref{e.dxgreen}.
\end{pf}

It is also important to note that
$\|\DX_{(D,E_1,E_2)}\| =\|\DX_{(D,\rev(E_2),\rev(E_1))}\|$; that is,
the total mass of $\DX_{(D,E_1,E_2)}$ is equal to that of
$\DX_{(D,\rev(E_2),\rev(E_1))}$. This is proved by reversing the paths, which
gives a measure-preserving bijection between the support
of these two measures.

%s6.2 ###
\subsection{Revisit probability estimate.}

We now return to the setup and notation of Section \ref{s.HEsetup}. Fix
some ball $B$ that intersects $\gamma[0,n]$. The present goal is to get
an upper bound for the conditional probability $\P[\gamma[n,N]\cap
B\ne\varnothing\md\gamma[0,n]]$, under some simple geometric
assumptions. To simplify notation, instead of discussing conditional
probabilities, we shall instead obtain a bound on $\P[\gamma[0,N]\cap
B\ne\varnothing]$ for a ball $B$ intersecting $\p D$. The conditional
probability estimate then readily follows, because of the Markovian
property mentioned in Remark~\ref{r.HEmark}.

\begin{proposition}\label{p.hitprob}
Let $0<r<R<\infty$. Let $B(z,r)$ be a ball of
radius $r$ intersecting $\p D$, and let $B(z,R)$ be
the concentric ball with radius $R$. Suppose that there is
no component of $B(z,R)\cap D$ whose
boundary intersects both $\AA_-$ and $\AA_+$.
Then
%
%e6.2 ###
\begin{equation}\label{e.hitprob}
\P\bigl[\gamma[0,N]\cap B(z,r)\ne\varnothing\bigr]
\le O(1) (r/R)^{\cc},
\end{equation}
where $\cc>0$ is a universal constant.
\end{proposition}

\begin{figure}

\includegraphics{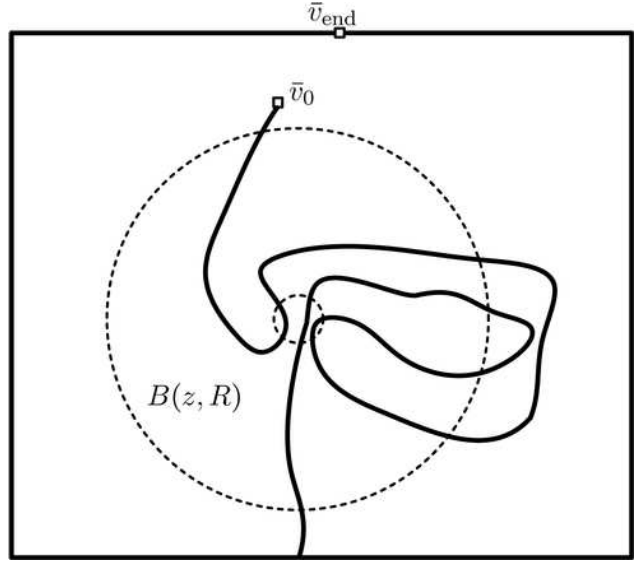}

\caption{A possible boundary of the domain $D$.}
\label{f.hitprob}
% , with $\oo_0$ and $\bo$ depicted as squares and $B(z,r)$ and
%$B(z,R)$ shown in dotted lines.}
\end{figure}

The reader should think about the case where $\p D$ is rather wild
geometrically, as in Figure~\ref{f.hitprob}. In particular, we need to
include the case where $B(z,r)$ intersects both $\AA_-$ and $\AA_+$.
The precise form of the right-hand side of~\eref {e.hitprob} will not
be important in the following. What is essential is only that it tends
to zero with $r/R$. If necessary, one can probably show that $\cc=1/2$
works, by using the following proof and a discrete version of the
Beurling projection theorem. (Some discrete version of the Beurling
projection theorem is given in~\cite{6}, Theorem~2.5.2, but it is not
precisely the same statement that would apply here.)

%fig3
%f3 ###

\begin{pf*}{Proof of Proposition~\ref{p.hitprob}}
We will assume that $r\ge5$. This involves no loss of generality,
because the distance from point on $\gamma[0,N]$ to $\p D$ is bounded from
zero, except near $\oo_0$ and $\bo$, and the result clearly holds
when $R$ is bounded.

Let $\ev Q$ be the event that there is a $j$ such that
$\oo_j\in B(z,r)$, and on $\ev Q$ let $\sigma$ be the
least such $j$. We assume, with no loss of generality,
that $R>9r$,
say, since otherwise the statement of the proposition is trivial.
Consider the ball $B(z,{3r})$.
By our assumptions, each component of $B(z,{3r})\cap D$
has boundary entirely in $\AA_-$ or entirely in $\AA_+$.
On the event $\ev Q$, let $S$ be the connected component
of $B(z,{3r})\cap D$ intersecting $\gamma[\sigma-1,\sigma]$.
Let $\ev Q_-\subset\ev Q$ be the event that $\p S\subset\AA_-$.
By symmetry, it is enough to prove that
$\P[\ev Q_-]\le O(1)(r/R)^{\cc}$.

Let $E_-$ denote the set of directed edges in $\inE=\inE(D)$
whose initial vertex is in $B(z,3r)\cap\HH_0^{-1}(0)$.
Let $D_n:=D\setminus V_n$, where $V_n$ is
as in Section~\ref{s.HEconstruction}
and let $E_+^n$ denote the set of directed edges connecting vertices
in ${D_n}$ to vertices in $\p D_n\cap\HH_n^{-1}(1)$.
The reason that the measures $\DX$ are useful here is because
the total mass $\|\DX_n\|$ of
$\DX_n:= \DX_{(D_n,E_-,E_+^n)}$ is a martingale.
One easy way to deduce this is by considering walks that hit
the vertex $v_n$ before any other vertex in $\p D_n$ (except for
the initial vertex of the walk).
Given the part of such a walk up to its first visit to $v_n$,
the probability that it first exits $D_{n-1}$ using an edge
from $E_+^{n-1}$ is precisely $h_{n-1}(v_n)$, which is just the
probability that $\HH_n(v_n)=1$.
Alternatively, Lemma~\ref{l.hmart} implies that
$\|\DX_n\|$ is a martingale,
because the total
mass of $\|\DX_n\|$ is just a linear combination of
the values of $h_n$ on the terminal vertices of $E_-$.

Since $\|\DX_n\|$ is a nonnegative martingale, the optional stopping
theorem implies that
%
%e6.3 ###
\begin{equation}\label{e.optional}
\|\DX_0\|\ge\E[ \mathbh{1}_{\ev Q_-}\|\DX
_\sigma\|].
\end{equation}
The bound on $\P[\ev Q_-]$ will be proven by estimating $\|\DX_0\|$ and
estimating $\|\DX_\sigma\|$
on the event $\ev Q_-$. By our assumption, every path starting from an
edge in $E_-$ and ending
with an edge in $E_+$ which stays in $D$ in between must exit $B(z,R)$.
Consequently,
$\|\DX_0\|\le\|\DX_{(D^R,E_-,E_R)}\|$, where $D^R$ is the intersection
of $D$ with $B(z,R)$
(adjusted to have as its boundary a subgraph of the grid $\mathit{TG}$), and
$E_R$ are the edges connecting
vertices in $D^R$ to vertices in $D\setminus D^R$. By Corollary~\ref
{c.hitball}, under the measure
$\DX_{(B(z,R),\rev(E_R))}$, the expected number of paths that hit
$B(z,R/2)$ is $O(1)$. For any
$r'\in(r,R)$, a random walk starting near the boundary of $B(z,r')$
has probability bounded away
from $0$ to exit $D$ before hitting $B(z,r'/2)$. Thus, the probability
that a random walk started
near the boundary of $B(z,R/2)$ will hit $B(z,3r)$ before exiting $D$
is $O(1)(r/R)^{\cc}$, for some
constant $\cc>0$. Hence,
%
%e6.4 ###
\begin{equation}\label{e.dx0}
\|\DX_0\|\le\bigl\|\DX_{(D^R,E_-,E_R)}\bigr\| =
\bigl\|\DX_{(D^R,\rev(E_R),\rev(E_-))}\bigr\|
\le O(1) (r/R)^{\cc}.
\end{equation}

Consider now the case where $\ev Q_-$ happens. As $\gamma[0,\sigma]$
crosses the annulus $B(z,{3r})\setminus B(z,r)$, there must be an arc
among the connected components of $\p B(z,2r)\cap
D\setminus\gamma[0,\sigma]$ that connects the right side of
$\gamma[0,\sigma]$ (where $\HH_n$ takes the value $1$) to the boundary
of $S$ (where $\HH_n$ takes the value $0$). By considering vertices
along this arc, we can find a vertex $v$ close to $\p B(z,2r)$ from
which the ratio between $H_{D_\sigma}(v,\rev(E_-))$ and
$H_{D_\sigma}(v,E_+^\sigma)$ is bounded and bounded away from zero by
universal constants (because these quantities do not vary by more than
a constant factor when moving from a vertex to its neighbor), or else
there is an edge in $E_-\cap E_+^\sigma$. In the latter case, clearly
$\|\DX_{(D_\sigma ,E_-,E_+^\sigma)}\|$ is bounded away
from zero. Consider therefore the case where such a $v$
exists. Since random walk starting from $v$ has probability bounded
away from zero to complete a loop going around the annulus
$B(z,3r)\setminus B(z,r)$ before exiting it, we have
$H_{D_\sigma}(v,\rev(E_-))+ H_{D_\sigma}(v,E_+^\sigma)$ bounded away
from zero. Consequently, each of these summands is bounded away from
zero. Let $B$ be the ball centered at $v$ whose radius is half the
distance from $v$ to $\p D_\sigma$. By Corollary~\ref{c.hitball}, the
measure under $\DX_{(D_\sigma,E_-)}$ of the set of paths hitting $B$ is
bounded away from zero. Since $H_{D_\sigma }(v,E_+^\sigma)$ is bounded
away from zero, a~random walk started at any vertex in $B$ has
probability bounded away from zero to exit $D_\sigma$ in $E_+^\sigma$,
by the Harnack principle (e.g., $k=0$ in \cite{9}, Lemma 5.2).
Consequently, we see that also in this case $\|\DX
_{(D_\sigma,E_-,E_+^\sigma)}\|$ is bounded away from zero on the event
$\ev Q_-$ by an absolute constant. Combining this
with~\eref{e.optional} and~\eref{e.dx0} establishes $\P[\ev
Q_-]=O(r/R)^{\cc}$. The proof for the event $\ev Q\setminus\ev Q_-$ is
entirely symmetric. \rightqed\end{pf*}

%s6.3 ###
\subsection{Local uniform convergence.}

\begin{proposition}[(Local uniform convergence)]\label{p.locunif}
In the setting of Theorem~\up{\ref{t.heunifconv}},
for every fixed $T>0$, there is a coupling
of $\hg$ and $\tg$ so that
\[
\sup\{|\hg(t)-\tg(t)|\dvtx 0\le t\le T\}\to 0
\]
in probability as $\rr\to\infty$.
\end{proposition}

\begin{pf}
Consider a sequence of pairs $(D,\phi)$, with $\rr\to\infty$. Using
Lemma~\ref{l.parthaus}, we
couple $\tg$ and the sequence $\hg$ so that for each $t\ge0$ the set
$\tg[0,t]$ is a.s. the
Hausdorff limit of $\hg[0,t]$.

Our strategy will be to prove that the curves $\hg[t, t^{\prime\prime}]$ converge to
$\tg[t, t^{\prime\prime}]$ in the
Hausdorff sense for all rational pairs $0 \leq t, t^{\prime\prime} \leq T$ (for
notational convenience, we may
assume $T$ is also rational). We will prove this, in turn, by showing
that for all rationals $t_0
< t <t'< t^{\prime\prime} < t_1 < T$, the Hausdorff limit of
$\hg[t,t^{\prime\prime}]$ is a.s.
disjoint from both $\tg[0,
t_0]$ and $\tg[t_1,T]$. We will begin by restricting our attention to a
large compact set and using
Proposition~\ref{p.hitprob} to derive upper bounds for the probability
that $\hg[s(s(t)),\infty)$
comes close to $\hg(s(t))$, where $t<s(t)<s(s(t))$ are chosen below so
that $\hg(s(t))$ is
``well exposed'' and so that the assumptions of
Proposition~\ref{p.hitprob} apply.

Let $\eps>0$, and let $F$ be some compact subset of $\H$. For $t\in
[0,\infty)$, let $s(t)$ be the
first $s>t$ such that $\hg(s)$ is in the unbounded connected component of
$\{z\in\H\dvtx \dist(z,\hg[0,t])\ge\eps\}$.
Fix some $t,\delta>0$,
and let $\ev A=\ev
A(t,\delta)$ be the event that $B(\hg(s(t)),2\eps)
\subset F$ and $\hg
[s(s(t)),\infty)\cap
B(\hg(s(t)),\delta)\ne\varnothing$. We claim that if $\rr$ is
sufficiently large, then
%
%e6.5 ###
\begin{equation}\label{e.aprob}
\P[\ev A]\le C_F (\delta/\eps)^{\cc},
\end{equation}
where $\cc$ is the same constant as in Proposition~\ref{p.hitprob},
and $C_F$ is a constant depending only on $F$.

By the Koebe distortion theorem, there is a constant $c=c(F)>0$
such that for every $z,z'\in F$ with $z\ne z'$ we have
%
%e6.6 ###
\begin{equation}\label{e.metriccompare}
c^{-1} \rr
\le
\frac{|\phi^{-1}(z)-\phi^{-1}(z')|}{|z-z'|}\le c \rr.
\end{equation}
[In other words, the metric in $F$ is comparable to the metric in
$\phi^{-1}(F)\subset D$ scaled by $\rr^{-1}$.] Let $n$ be the least
integer such that $\phi(\gamma[0,n])\supset\hg[0,s(s(t))]$. (Recall
that the parameterization of $\gamma$ is not by capacity, but is
proportional to arclength.) Now condition on $\gamma[0,n]$, and assume
that $B(\hg(s(t)),2\eps)\subset F$. (If~that has zero probability, then
$\P[\ev A]=0$.) Let $z_0:=\phi^{-1}(\hg(s(t)))$. Assume that $\rr$ is
sufficiently large so that $|\phi(\gamma(n))-\hg(s(s(t)))|<\eps/3$;
indeed, how large $\rr$ is required to be can be determined from the
constant $c$ in~\eref{e.metriccompare}. The metric
comparison~\eref{e.metriccompare} implies that
$\phi(B(z_0,\rr\eps/(2c)))\subset B(\hg(s(t)),\eps/2)$. This and the
\mbox{definition} of $s(s(t))$ imply that there is a path in
$D\setminus(\gamma[0,n]\cup B(z_0,\rr\eps/(2c)))$ connecting
$\gamma(n)$ and $\bo$. We may then apply Proposition~\ref{p.hitprob}
with $D$ replaced by the $D_n$ of Remark~\ref{r.HEmark} to conclude
that conditioned on $\gamma[0,n]$ the probability that $\gamma[n,N]$
intersects $B(z_0,\delta'\rr)$ is $O(1)(2c\delta'/\eps)^{\cc}$.
Now~\eref{e.aprob} follows by another application
of~\eref{e.metriccompare}.

Let $\Gamma^{t'}_t$ denote the Hausdorff limit of $\hg[t,t']$, when it exists.
Note that we may pass to a subsequence of pairs
$(D,\phi)$ so that $\hg
[t,t']$ converges in law.
Consequently, by using
\cite{3}, Theorem 11.7.2  again, we may
assume that almost
surely the limits $\Gamma^{t'}_t$,
$\overline s_t:=\lim_{(D,\phi)} s(t)$,
$\overline{ss}_t:= \lim_{(D,\phi)} s(s(t))$ and $p_t:=\lim_{(D,\phi)}
\hg(s(t))$ exist for every
pair of rationals $0<t<t'<\infty$.

Let $t>0$ be rational, and let $S(t,b)$ be the set of $t'>t$ such that
$\tg(t')$ is in the unbounded component of
$\{z\in\H\dvtx \dist(z,\tg[0,t])\ge b\}$.
By construction, it is clear that
$\dist(p_t,\tg[0,t])=\eps$,
$p_t\in\tg[0,\overline{ss}_t]$,
$\overline s_t< \inf S(t,2\eps)$
and $\overline{ss}_t<\inf S(t,3\eps)$.
Let $t'$ be a rational satisfying
$\overline{ss}_t<t'<\inf S(t,3\eps)$.

Now, the following a.s. statements will hold on the event $\inf
S(t,3\eps)<T$ and $B(p_t,3\eps)\subset F$. First, by \eref{e.aprob}, we
have a.s. $p_t\notin\Gamma^{T}_{t'}$. Note that
$p_t\in\tg[0,t']\setminus\tg[0,t]$. Since $\tg$ is a simple path,
$\tg[0,T]\setminus p_t$ has two disjoint components, one containing
$\tg[0,t]$, and the other containing $\tg(T)$. Since $\Gamma^{T}_{t'}$
is connected and contains $\tg(T)$, we conclude that a.s.
$\Gamma^{T}_{t'}\cap\tg[0,t]=\varnothing$ [on the event $B(p_t,3\eps
)\subset F$, $\inf S(t,3\eps)<T$]. Now let $t^{\prime\prime}$ be any
rational in $(t,T)$. Since $\tg$ is a simple path, there is some small
$\eps>0$ such that $t^{\prime\prime}\in S(t,3\eps)$.

Because $F$ is an arbitrary compact subset of $\H$, $p_t\in\H$ a.s.
and in the above discussion
$\eps>0$ was arbitrary, it follows that
$\Gamma^{T}_{t^{\prime\prime}}\cap\tg
[0,t]=\varnothing$ a.s. Since
$\Gamma^{t^{\prime\prime}}_t\subset\tg[0,T]$ and is a.s.
disjoint from $\tg
[0,t_0]\cup\tg[t_1,T]$ whenever
$t_0<t<t^{\prime\prime}<t_1$ (by the above), it follows that
$\Gamma^{t^{\prime\prime}}_t=\tg
[t,t^{\prime\prime}]$ a.s. for every pair of
rationals $0<t<t^{\prime\prime}<T$. Consequently,
$\Gamma^{t^{\prime\prime}}_t=\tg[t,t^{\prime\prime}]$ a.s.
for every pair $0\le t\le
t^{\prime\prime}\le T$
(the cases $t=0$ and $t^{\prime\prime}=T$ are similarly treated). Thus,
$\lim_{(D,\phi)}
\sup_{t\in[0,T]} |\hg(t)-\tg(t)|=0$. Since every sequence of $(D,\phi)$
with $\rr\to\infty$ has a
subsequence such that this holds, this also holds without passing to a
subsequence. The proposition
is thus established.
\end{pf}

%s6.4 ###
\subsection{Uniform transience and conclusion.}

Since $\phi(\bo)=\infty$, $\hg$ is transient; that is, $\lim_{t\to\infty
}|\hg(t)|=\infty$. The
following is a uniform version of this statement.

\begin{proposition}\label{p.uniftrans} For every $\eps>0$ and $R>0$
there is a $T=T(R,\eps)$ such
that
\[
\P\bigl[\hg[T,\infty)\cap B(0,R)\ne\varnothing\bigr]<\eps
\]
if $\rr$ is sufficiently large.
\end{proposition}

The reader may note the clear similarity with Proposition~\ref{p.hitprob}.
The main difference is that there the path considered was in the domain $D$,
whereas here the path $\hg$ is in the image under the conformal map,
that is,
in the upper half plane $\H$. Indeed, the proof is quite similar.

The following lemma about the excursion measures will be needed.

\begin{lemma}\label{l.phiDX}
Let $R>0$. Suppose that at each vertex $v\in
V(\mathit{TG})\cap\phi^{-1} (\overline\H\cap\overline{B(0,R)})$ a
simple random walk $X^v$ is started, and the walk is stopped at the
first time $t\ge1$ such that $X^v(t)$ exits $D$ or
$|\phi(X^v(t))|\notin(R,2R)$. Then the expected number of walks which
stop when $|\phi(X^v(t))|\ge2R$ is $O(1)$.
\end{lemma}

\begin{pf}
Let $V_D:=V(\mathit{TG})\cap\closure{D}$. For $v\in V_D$, let $Y^v$
denote a simple random walk on $V_D$, where at each step the walk jumps
with equal probability along each of the edges $e$ with
$e\subset\closure D$ and $Y^v(t)\in e$. Let $V_0:=\{v\in V_D\dvtx
|\phi(v)|\le R\}$ and $V_1:=\{v\in V_D\dvtx |\phi(v)|\ge2 R\}$. Let $M$
be the expected number of walks $Y^v$ with $v\in V_0$ such that
$Y^v(\sigma)\in V_1$, where $\sigma:=\inf\{t\ge1\dvtx Y^v(t)\in V_0\cup
V_1\}$. It clearly suffices to show that $M=O(1)$. (The difference from
the $X^v$ is that the $Y^v$ are reflected off of $\p D$, rather than
killed there.)

For a function $f$ on $V_D$ and $v\in V_D$, let
$\Delta f(v):=\sum_{[v,u]} (f(u)-f(v))$,
where the sum extends over edges $e=[v,u]$ containing $v$ such that
$e\subset\closure D$.
Now let $g\dvtx V_D\to[0,1]$ be the unique function such
that $\Delta g(v)=0$ for $v\in V_D\setminus(V_1\cup V_2)$,
$g=0$ on $V_0$ and $g=1$ on $V_1$; that is, $g$ is harmonic
in $V_D\setminus(V_0\cup V_1)$ with
the appropriate boundary values. It is immediate to verify
that
\[
M=\sum_{v\in V_0}\Delta g(v)/d_v\le\sum_{v\in V_0} \Delta g(v),
\]
where
$d_v$ is the number of edges containing $v$ in $\closure D$,
since $g(v)$ is
the probability that $Y^v$ hits $V_1$ before $V_0$.
Because $\sum_{v\in V_D} \Delta g(v)=0$,
we have
\[
M\le-\sum_{v\in V_1}\Delta g(v)\,.
\]

Let $E(f):=\sum_{[v,u]} (f(v)-f(u))^2$,
where the sum runs over all edges in $\closure D$.
It is well known (and simple to show) that $g$ minimizes
$E$ among functions mapping~$V_1$ to $1$ and $V_0$ to $0$.
For each edge $[v,u]$, we may distribute the
quantity $(f(v)-f(u))^2$ by giving $f(v)(f(v)-f(u))$
to the vertex $v$ and $f(u)(f(u)-f(v))$ to the vertex $u$.
Consequently, by summing over the contributions to each vertex, we find
\[
E(f)=-\sum_v f(v)\Delta f(v).
\]
Hence, $E(g)=-\sum_{v\in V_1} \Delta g(v)$, which gives
$M\le E(g)$.

There are several different ways to estimate $E(g)$ and complete the proof.
We opt for an easy and short argument,
which unfortunately does require terminology and
results from the literature.
(Similar arguments appear, e.g.,
in~\cite{2,5,11}.)

Let $A$ be the annulus $A:=B(0,2 R)\setminus B(0,R)$. Let $\tilde g(v)
:= (|\phi(v)|-R)/R$ when $v\in V_D\setminus(V_0\cup V_1)$ and $\tilde
g(v) := j$ when $v\in V_j$, $j=0,1$. Then $E(g)\le E(\tilde g)$, by the
characterization of $g$ as a minimizer. Following~\cite{14}, we say
that a set $F\subset\C$ is $s$-fat if for every disk $B=B(z,r)$ with
$z\in F$ and $F\not\subset B$, we have $\area(F\cap B)\ge s\area(B)$.
Consider a triangle $\triangle$ of the grid $\mathit{TG}$ with
$\triangle\subset\closure D$. Then $\phi(\triangle)$ is a
$K$-quasidisk, for some constant $K$, by~\cite{4}. [There, it is
required that $\triangle\subset D$, but $\triangle\subset \closure D$
works too, by standard compactness properties. Besides, it suffices for
the argument given below that $\phi(\triangle')$ is a $K$-quasidisk
when $\triangle'$ is a slightly rescaled copy of $\triangle$ contained
in the interior of $\triangle$, provided that $K$ does not depend on
the scaling factor.) Consequently, by~\cite{14},  Corollary~2.3,
$\phi(\triangle)$ is $s$-fat for some constant $s>0$. Thus,
$\diam(\phi(\triangle)\cap A)^2\le O(1)\area(\phi(\triangle )\cap B(0,3
R))$. Hence, $\sum_\triangle\diam(\phi(\triangle)\cap A)^2 = O(R^2)$,
where the sum extends over all the triangles of the grid $\mathit{TG}$
that are contained in $D$. Now, if $[u,v]$ is an edge in $\closure D$,
then there is a grid triangle $\triangle\subset\closure D$ with
$[u,v]\subset\p\triangle$. We have then $(\tilde g(v)-\tilde
g(u))^2\le\diam(\triangle\cap A)^2/R^2$. Since each triangle has three
edges, we get
\[
M\le E(g)\le E(\tilde g) \le3
\sum_\triangle\diam(\triangle\cap
A)^2/R^2=O(1).
\]
This completes the proof.
\end{pf}

\begin{pf*}{Proof of Proposition~\ref{p.uniftrans}}
We choose $R_1=R_1(\eps,R)$ much larger than~$R$.
Let $t_1:=\inf\{t\ge 0\dvtx |\hg(t)|=R_1\}$, and let
$T$ be a constant such that
\[
\P[T-1<\inf\{t\ge0\dvtx |\tg(t)|=R_1\}]<\eps/3.
\]
Then, by Lemma~\ref{l.parthaus}, when $\rr$ is sufficiently large, we have
$\P[T<t_1]<\eps/3$.
Consequently, it suffices to show that
%
%e6.7 ###
\begin{equation}\label{e.tworads}
\P\bigl[\hg[t_1,\infty)\cap B(0,R)\ne\varnothing\bigr]<\eps/2.
\end{equation}
Let $m$ be the least integer such that
$\phi\circ\gamma[0,m]\supset\hg[0,t_1]$.
%We now condition on $\gamma[0,m]$.
The proof now proceeds as in Proposition~\ref{p.hitprob}, with only
minor changes, which will be henceforth described. Let $\ev Q$ be the
event that there is a $j>m$ such that $|\phi(\gamma(j))|<R$, and on
$\ev Q$ let $\sigma$ be the least such $j$. As in the proof of
Proposition~\ref{p.hitprob}, the event $\ev Q_-\subset\ev Q$ is
defined. For integer $n\in[m,N]$, we consider the excursion measure in
$D\setminus V_n$ with excursions started at the vertices in
$\alpha:=\HH_m^{-1}(0)\cap\phi^{-1}(B(0,3R))$ and terminating at
vertices in $\beta_n:=\HH_n^{-1}(1)$. The total mass of this measure is
a martingale. It suffices to show that this is very small at $n=m$, but
is bounded away from zero at $n=\sigma$ on the event $\ev Q_-$.

We first do the estimate for $n=m$.
The expected number of excursions in $D$ from $\alpha$ that hit
$\phi^{-1}(\p B(0,6R))$ is the same as the number
of excursions in the domain which consists of
the grid triangles intersecting $\phi^{-1}(B(0,6R))$ starting at
vertices in
$\phi^{-1}(\closure\H\setminus B(0,6R))$ that hit $\alpha$,
by symmetry, and this quantity is bounded by
the number of excursions in the domain which is
essentially $\phi^{-1}(B(0,6R)\setminus B(0,3R))$ starting
at vertices in
$\phi^{-1}(\closure\H\setminus B(0,6R))$
that hit $\phi^{-1}(B(0,3R))$.
By symmetry again, this is the same as the expected number of
excursions in the reverse direction. This quantity is $O(1)$,
by Lemma~\ref{l.phiDX}.
Consequently, the expected number of excursions from $\alpha$
in $D$ that cross $\phi^{-1}(\p B(0,6R))$ is $O(1)$.
It is not hard to see that when $\rr$ is sufficiently large,
there will not be any grid edge crossing both
$\phi^{-1}(\p B(0,6R))$ and $\phi^{-1}(\p B(0,7R))$,
for example, by considering the harmonic measure from~$\bb$ of such
an edge.
Now, \cite{9}, Lemma 5.4 tells us that a random walk
started in $\phi^{-1}(B(0,7R))$ has probability $o(1)$
to exit $\phi^{-1}(B(0,R_1))$ before exiting $D$, uniformly as
$R_1/R\to\infty$. (That lemma refers to the square grid, but the
proof applies here as well. Also, in that lemma the image conformal
map is onto the unit disk $\U$, but this is simply handled by choosing an
appropriate conformal homeomorphism from $\U$ to $\H$.)

It remains to prove a bound from below for measure of excursions in
$D\setminus V_\sigma$ from $\alpha$ to $\beta_\sigma$, on the event
$\ev Q_-$. As in the proof of Proposition~\ref{p.hitprob}, it suffices
to find a vertex $v$ such that the discrete harmonic measure in
$D\setminus\gamma[0,\sigma]$ from $v$ of each of the sets $\alpha$ and
the left side of $\gamma[0,\sigma]$ is bounded from below. Consider any
vertex $w$ near $\phi^{-1}(\p B(0,3R))$. The continuous harmonic
measure from~$\phi(w)$ of~$\R$ in the domain $\H \cap(B(0,4R)\setminus
B(0,2R))$ is bounded from below. By the convergence of discrete
harmonic measure to continuous harmonic measure, when $\rr$ is large, a
random walk started at $w$ will have probability bounded from below to
hit $\p D$ before exiting $\phi^{-1} (B(0,4R)\setminus B(0,2R))$.
(Specifically, while not close to the boundary, the random walk behaves
like Brownian motion, which is conformally invariant. Once it does get
close to the boundary, we may apply~\cite{9}, Lemma~5.4, say.) As in
the proof of Proposition~\ref{p.hitprob}, on $\ev Q_-$ we can find a
vertex~$v$ near $\phi^{-1}(\p B(0,3R))$ where the discrete harmonic
measure of $\alpha$ is comparable to that of $\beta_\sigma$. Hence,
both are bounded away from zero. This completes the proof.
\rightqed\end{pf*}

\begin{pf*}{Proof of Theorem~\ref{t.heunifconv}}
The theorem follows immediately from
Propositions \ref{p.locunif}~and~\ref{p.uniftrans}.
\end{pf*}

\section*{Acknowledgments.}

We wish to thank Richard
Kenyon and David Wilson for
inspiring and useful conversations.

\printaddresses

\end{document}